\documentclass{amsproc}

\newtheorem{theorem}{Theorem}[section]
\newtheorem{lemma}[theorem]{Lemma}
\newtheorem{proposition}[theorem]{Proposition}
\newtheorem{corollary}[theorem]{Corollary}

\theoremstyle{definition}
\newtheorem{definition}[theorem]{Definition}

\theoremstyle{remark}
\newtheorem{remark}[theorem]{Remark}

\numberwithin{equation}{section}

\newcommand{\R}{\mathbb R}

\newcommand{\N}{\mathbb N}
\newcommand{\C}{\mathbb C}
\newcommand{\del}{\partial}

\newcommand{\ds}{\displaystyle}
\newcommand{\Par}{\mathcal P}
\newcommand{\calC}{\mathcal C}
\newcommand{\calL}{\mathcal L}

\copyrightinfo{2002}
  {American Mathematical Society}

              [1996/10/25 v1.2b CONM-P Author Class]


\begin{document} 

\title{Bifurcating nodoids}
\author{Rafe Mazzeo}
\address{Stanford University}
\email{mazzeo@math.stanford.edu}
\thanks{Supported by the NSF under grant DMS-9971975 and at 
MSRI by NSF grant DMS-9701755} 

\author{Frank Pacard}
\address{Universit\'e Paris XII}
\email{pacard@univ-paris12.fr}
\thanks{Supported at MSRI by NSF grant DMS-9701755}

\begin{abstract}
All complete, axially symmetric surfaces of constant mean curvature
in ${\R^3}$ lie in the one-parameter family $D_\tau$ of Delaunay surfaces. 
The elements of this family which are embedded are called unduloids;
all other elements, which correspond to parameter value $\tau \in 
{\mathbb R}^-$, are immersed and are called nodoids. The unduloids 
are stable in the sense that the only global constant mean curvature 
deformations of them are to other elements of this Delaunay family. We 
prove here that this same property is true for nodoids only when $\tau$ 
is sufficiently close to zero (this corresponds to these surfaces  
having small `necksizes'). On the other hand, we show that as $\tau$ decreases 
to $-\infty$, infinitely many new families of complete, cylindrically bounded 
constant mean curvature surfaces bifurcate from this Delaunay family. 
The surfaces in these branches have only a discrete symmetry group. 
\end{abstract}

\maketitle 

\section{Introduction}
In 1841, C. Delaunay discovered a beautiful one-parameter family of 
complete noncompact surfaces of constant mean curvature one in 
$\R^3$ which are invariant under rotations about an axis \cite{Del}. 
(Henceforth we shall abbreviate 
constant mean curvature one by CMC.) Using the rotational symmetry,
the search for these surfaces reduces to finding their meridian curves,
and these may be found, in turn, by solving an appropriate ODE. 
This family of `Delaunay surfaces' is parametrized by a variable
$\tau$ lying in the set
\[
\Par = (-\infty,1] - \{0\} = \Par^+ \cup \Par^-; \qquad 
\Par^- = (-\infty,0), \qquad \Par^+ = (0,1].
\]
The Delaunay surface corresponding to any value $\tau \in \Par$
will be denoted $D_\tau$. When $\tau \in\Par^+$, then $D_\tau$ 
is embedded and is called an unduloid; the meridian curve in this
case is the roulette of an ellipse. When $\tau \in \Par^-$, then
$D_\tau$ is no longer embedded, and is called a nodoid; the meridian
curve of any nodoid is the roulette of a hyperbola. There is 
a good geometric limit of the surfaces $D_\tau$ as
$\tau \to 0$ (which is the same for $\tau \searrow 0$ and 
$\tau \nearrow 0$), which is the nodded surface formed by the infinite
union of mutually tangent spheres of radius $1$ arranged along
a common axis. A nice geometric description of these Delaunay surfaces
can be found in \cite{Eel}.  

In this paper we shall prove the existence of some rather surprising 
new families of complete CMC surfaces which arise as deformations
of Delaunay nodoids. Let us digress briefly to set these new surfaces 
into their proper context before describing our results in more detail.

The past several years have witnessed great progress in the understanding 
of more general complete, CMC surfaces of finite topology. On the 
one hand, the great variety and flexibility of such surfaces has been
suggested by numerous numerical and computer experiments, notably
by the GANG group at the University of Massachusetts and Grosse-Brauckmann
and his collaborators at the Universit\"at Bonn, see
http://www.gang.umass.edu/cmc/ for many nice pictures. On the
other hand, there have been several theoretical advances which
vindicate many of these computer experiments. These advances
may be subdivided into two groups of results, those concerning
various constructions of complete CMC surfaces, which are not relevant
to the present discussion, and those concerning the structure theory 
of these surfaces.

This structure theory
is much better developed when we restrict attention to complete 
(oriented) Alexandrov-embedded CMC surfaces of finite topology, the 
key point being that any one of these surfaces has only finitely many
ends, and each of these ends is necessarily modeled on a Delaunay 
unduloid. Recall first that any surface of this type 
is conformally equivalent to a punctured Riemann surface $\Sigma = 
\overline{\Sigma}\setminus \{p_1, \ldots, p_k\}$; let $Y^3$ be a
handlebody such that $\del Y = \overline{\Sigma}$. Then $\Sigma$ 
is said to be Alexandrov embedded if its immersion into $\R^3$
extends to an immersion of $Y\setminus \{p_1,\ldots, p_k\}$. 
In particular, $D_\tau$ is not Alexandrov-embedded when $\tau < 0$. 

A deep theorem of Meeks \cite{Mee} states that any end of an 
Alexandrov-embedded CMC surface of finite topology is cylindrically 
bounded. Following this, Korevaar, Kusner and Solomon \cite{Kor-Kus-Sol} 
showed that for any such end there is a Delaunay unduloid $D_\tau$ 
to which the end converges exponentially. \cite{KK} contains a 
further global structure theorems. These results show that
Delaunay unduloids play a fundamental r\^ole as building blocks
for finite topology Alexandrov-embedded CMC surfaces. In addition, 
this tameness of the geometry of the ends leads to a fairly detailed 
understanding of many analytic and geometric problems on these 
surfaces, \cite{Kus-Maz-Pol}, \cite{Maz-Pac-Pol-2}. 

It is natural to try to extend this theory to complete, CMC surfaces 
with finite topology which are not necessarily Alexandrov-embedded,
for example those with ends modeled on Delaunay nodoids. Fairly 
general examples of such surfaces are constructed in \cite{Maz-Pac}, 
\cite{Maz-Pac-Pol-2}, but the ends of those surfaces are modeled on 
Delaunay nodoids with $\tau$ very close to zero. 

The results of this paper begin to elucidate the limitations of the 
structure theory \cite{Kor-Kus-Sol}, and how it fails when 
Alexandrov-embeddedness is dropped. More specifically, we examine
the stability theory of the Delaunay nodoids $D_\tau$ as $\tau \to
-\infty$. In a sense we explain below, there is some stable range
$\tau_* < \tau < 0$ for which the only deformations of $D_\tau$
are to other Delaunay surfaces. However, as $\tau$ decreases,
we establish the existence of new CMC surfaces which bifurcate
from the Delaunay family. These surfaces have only a discrete
symmetry group, rather than full rotational symmetry, and the
number and geometric complexity of these bifurcations increases as
$\tau$ becomes more negative. All of these new surfaces remain
cylindrically bounded, but their geometric structure is rather
intricate. 

To conclude this general discussion, we recall that the basic
theory of the moduli space of complete Alexandrov-embedded
CMC surfaces of finite topology is established in \cite{Kus-Maz-Pol}.
The paper \cite{Maz-Pac-Pol-2} extends this theory to include CMC surfaces with
ends modeled on Delaunay nodoids with $\tau > \tau_*$; (it also
develops some global aspects of this moduli space theory).
Our results here show that a comprehensive description of the 
moduli space theory of complete immersed CMC surfaces of finite topology
must incorporate these new families of CMC surfaces. 

\medskip

We may now give the precise statements of our main results. These
require a few preliminary definitions.
\begin{definition} We set terminology for a family of rigid
motions in $\R^3$, normalized for convenience to fix the $z$-axis, 
as well as the surfaces which have the associated symmetries. 
\begin{itemize}
\item[i)] Let $j \in \N, \ j \geq 2$. Let $R_j$ denote the rotation by 
angle $2\pi/j$ about the $z$-axis in $\R^3$. We say that a surface 
$\Sigma\subset {\R}^3$ is $R_j$-symmetric if it is invariant by $R_j$, but 
not by any $R_{j'}$ for $j' > j$.
\item[ii)] For $\alpha \in [- \pi, \pi]$, let $R_\alpha$ denote the rotation 
by angle $\alpha$ in $\R^2$. Also, let $t_\alpha \in \R$. We say that
the rigid motion $S_\alpha$ of $\R^3$ is a screw motion about the $z$-axis
(of angle $\alpha$ and translation length $t_\alpha$) if it has the form
\[
\begin{array}{rcccc}
S_\alpha & : &  \R^2\times \R & \longrightarrow & \R^2 \times \R \\
              &  &   (x_1,x_2,z)    & \longmapsto  & (R_\alpha\, x, z+t_\alpha)
\end{array}
\]
We say that the surface $\Sigma \subset \R^3$ is $S_\alpha$-symmetric 
if it is invariant with respect to some screw motion $S_\alpha$, for
some choice of translation length $t_\alpha$, and if it is not 
invariant with respect to any other screw motion $S_{\alpha'}$ with
$0 < \alpha' < \alpha$ if $\alpha \in (0,\pi]$ or $\alpha < \alpha' < 0$
if $\alpha \in [-\pi,0)$. 
\item[iii)] Finally, we say that the surface $\Sigma \subset \R^3$ is 
$T_{j,\alpha}$-symmetric if it is both $R_j$-symmetric and
$S_\alpha$-symmetric.
\end{itemize}
\label{de:syms}
\end{definition}

Our first main result states that there are infinitely many bifurcations
from the Delaunay family $D_\tau$ as $\tau \to -\infty$ to new 
CMC surfaces which have a discrete symmetry group:
\begin{theorem}
Let $j\geq 2$; then there for any $\alpha 
\in [-\pi/j,\pi/j]- \{0\}$, there exists a $T_{j,\alpha}$-symmetric 
CMC surface. In fact, each of these surfaces lies in a branch which
bifurcates from the nodoid  $D_{\tau_j}$ at some value $\tau_j <-\sqrt{j^2-2}$.
\label{th:2.1}
\end{theorem}
This may be proved either using the general bifurcation theorem of Smoller and 
Wasserman \cite{Smo-Was} or more simply by degree theory \cite{Smo}. As in most
such problems, the main point is to show that the index (i.e.\ number of 
negative eigenvalues) of the Jacobi operator $L_\tau$ on some compact
quotient of $D_\tau$ increases as $\tau \searrow -\infty$.  
As part of this, we show that there exists some value $\tau_* \in 
(-2,-\sqrt 2)$ such that no bifurcations occur when $\tau_* < \tau < 0$. 

\medskip

This result can be sharpened when $j$ is large enough, for in this case
we can guarantee that the bifurcation value $\tau_j$ is large negative.
This allows us to obtain better control on the spectrum of Jacobi 
operator, and so we may use the bifurcation theorem of Crandall and 
Rabinowitz \cite{Cra-Rab} to show that the bifurcation branches are smooth. 
\begin{theorem}
There exists $j_0 \geq 2$ such that when $j\geq j_0$ and $\alpha  \in 
[-\pi/j, \pi/j]- \{0\}$, then there are values $\tau_{j,\alpha} < 0$ 
$\eta_{j,\alpha} > 0$, and a real analytic branch of 
$T_{j,\alpha}$--symmetric CMC surfaces, which we denote $D_{j,\alpha}(
\tau,\eta)$, parametrized by $\eta$ with $|\eta| \leq \eta_{j, \alpha}$, 
which bifurcate from $D_{\tau_{j,\alpha}}$. These families of surfaces have
the following properties:
\begin{enumerate}
\item[(i)] As $\eta \to 0$, the surfaces $D_{j,\alpha}(\eta)$ converge 
uniformly on compact subsets to the nodoid 
$D_{\tau_{j,\alpha}}$, in $\calC^\infty$ topology. 
\item[(ii)] Locally on $D_{\tau_{j, \alpha}}$ and for $\eta$ sufficiently
small, we can write $D_{j,\alpha}(\eta)$ as a normal graph over 
$D_{\tau_{j,\alpha}}$ for some function $w_\eta$. The functions 
$w_\eta /\eta$ converge uniformly on compacts sets to a 
nontrivial function of the form $(s, \theta)  \rightarrow 
\phi (s)\,\cos (j\theta)$. (The variables $(s,\theta)$ will be defined in 
\S 2.) Furthermore, the 
function $\psi(s,\theta) := \phi (s)\, e^{i\theta}$ is a nontrivial 
solution of $\calL_\tau \psi =0$ which satisfies 
$\phi(s+2 \,\pi\,s_\tau,\theta) = e^{i\,j \, \alpha}\,\phi(s,\theta)$.
\end{enumerate}
\label{th:2.2}
\end{theorem} 
In fact, we prove that
\begin{equation}
\tau_{j, \alpha} = -\frac{j}{\sqrt{1 - (\alpha j/2\pi)^2}}
+ {\mathcal O}(j^{-1}).
\label{eq:2.1}
\end{equation}
\begin{remark}
It is likely that Theorem~\ref{th:2.2} holds whenever $j \geq2$. 
We also suspect that for any $j \geq 2$, $|\alpha| \leq \pi/j$,
$\alpha \neq 0$, there are precisely two bifurcating branches,
and these correspond to two separate bands of continuous spectrum
of $\calL_\tau$, with one branch bifurcating at the value 
$\tau = \tau_{j, \alpha}$, which lies in the first band, and 
the other bifurcating at a different value $\bar \tau_{j, \tau}$ lying 
in the second band. This second bifurcation point should satisfy
\[
\bar\tau_{j,\alpha} = 
-\frac{j}{\sqrt{1 - (|\alpha j|/2\pi - 1)^2}}
+ {\mathcal O}(j^{-1}).
\]
as $j \to \infty$.
\end{remark}

We conclude this introduction by mentioning that the bifurcations
we find are related to the so-called `Rayleigh instability of
the cylinder'. We refer in particular to \cite{Chand}, \S 111, for a 
discussion of this phenomenon as manifested in the capillary instability 
of a liquid jet. As an historical aside, this problem was originally 
studied by Plateau, who posited that a cylindrical jet should break 
up into rotationally symmetric pieces. Lord Rayleigh found fault with 
Plateau's argument but emended it by noting that in certain situations 
the instability should produce nonrotationally symmetric perturbations.

Other types of `CMC cylinders' with few or no symmetries are known
to exist. For example, it is pointed out in \cite{Kor-Kus-Sol} that 
a `relaxation' of the construction of Wente tori produces CMC surfaces 
with finite rotational and discrete translational symmetry, or with 
discrete screw motion symmetry. In addition, one can produce immersed 
CMC surfaces with continuous screw motion symmetry using the classical 
associate family construction. Finally, \cite{KMS} uses the 
DPW (Dorfmeister-Pedit-Wu) method (which is akin to the Weierstrass 
representation formula for minimal surfaces) to produce examples of 
cylindrically bounded CMC surfaces with no symmetry. To our knowledge,
none of these lie in a continuous family including the Delaunay surfaces. 

We are grateful to the referee for drawing to our attention 
the work of Lord Rayleigh, as well as the examples in the last paragraph.

\section{Isothermal parametrization of Delaunay surfaces}

Since the Delaunay surfaces $D_\tau$ are surfaces of revolution, their
most natural parametrizations would seem to be the obvious ones:
\[
(u,\theta) \longmapsto (\rho(u)\cos\theta,\rho(u)\sin\theta,u),
\]
where function $\rho = \rho_\tau$ is a solution of an ODE which
is derived from the constant mean curvature condition. However, for
most analytic purposes it turns out to be far more convenient to
use a different parametrization which is isothermal, and which we
now describe. This definition may seem ad hoc, but is motivated by a 
systematic line of reasoning in the theory of integrable systems; 
\cite{Maz-Pac} contains a detailed derivation of the fact 
that this parametrizes $D_\tau$.  

This isothermal parametrization rests on two functions 
$\sigma = \sigma_\tau$ and $\kappa = \kappa_\tau$, 
the definitions of which vary, according to whether $\tau$ is positive 
or negative. First, $\sigma$ is the unique smooth nonconstant solution to the initial 
value problem
\begin{align}
\left(\frac{d\sigma}{ds}\right)^2 + \tau^2 \cosh^2 \sigma & =1, \qquad 
\del_s\sigma(0) = 0, \quad \sigma (0) < 0, \quad \mbox{when}\ 
\tau \in \Par^+, \label{eq:1.2} \\
\left(\frac{d\sigma}{ds}\right)^2 + \tau^2 \sinh^2 \sigma & =1, \qquad 
\del_s\sigma(0) = 0, \quad \sigma (0) < 0, \quad \mbox{when}\ 
\tau \in \Par^- .
\label{eq:1.21}
\end{align}
Next, $\kappa(s)$ is the unique solution of 
\begin{align}
\frac{d\kappa}{ds} & = \ \tau^2 \, e^{\sigma}\, \cosh \sigma, \qquad 
\kappa(0) = 0,
\quad \mbox{when}\ \tau \in \Par^+, \label{eq:1.4} \\
\frac{d\kappa}{ds} \ &= -\tau^2\, e^{\sigma}\,\sinh\sigma, 
\qquad \kappa(0) = 0, \quad \mbox{when}\ \tau \in \Par^-. 
\label{eq:1.41}
\end{align}

The change of variables from the previous cylindrical coordinates
$(u,\theta)$ is effected by setting $u = \kappa(s)$ and one then
has (see \cite{Maz-Pac}) $\rho (\kappa(s)) = \tau e^{\sigma(s)}$, and so
the new isothermal parametrization is given by
\begin{equation}
X_\tau: \R \times S^1 \ni  (s,\theta) \longmapsto 
\frac12 \, \left(\tau\,e^{\sigma_\tau(s)} \, \cos\theta,\tau\,
e^{\sigma_\tau(s)}\,\sin\theta,\kappa_\tau (s)\right).
\label{eq:1.1}
\end{equation}
The metric coefficients in this new coordinate system
are $g_{ss} = g_{\theta \theta} = \tau^2 e^{2\sigma}$, 
$g_{s \theta}= g_{\theta s} = 0$. 

Regardless of the sign of $\tau$, the function $s \to \sigma_\tau(s)$ 
necessarily changes sign. Together with (\ref{eq:1.41}), this means that 
when $\tau < 0$, then $\kappa$ is not monotone; in contrast, using
(\ref{eq:1.4}), $\kappa$ is monotone when $\tau > 0$. 
Hence although the differential of $X_\tau$ is always full rank, so that 
$X_\tau$ is always an immersion, it is only an embedding when $\tau > 0$.  
Notice also that when $\tau = 1$, then (\ref{eq:1.2}) implies that 
$\sigma \equiv 0$, and so $\kappa(s) = s$ and $\rho  \equiv 1/2$, which 
means that $D_1$ is a cylinder of radius $1/2$. 

\section{The period function}

The Hamiltonian nature of the equations (\ref{eq:1.2}) and 
(\ref{eq:1.21}) implies that $\sigma $ is periodic. We denote its 
period by a
$2\pi s_\tau$. We now investigate the dependence of $s_\tau$ on $\tau$.

To begin, $s_\tau$ has an integral representation. To state it, 
let $a_\tau, a_\tau' > 0$ be determined by the equations
$\tau \cosh a_\tau = 1$ when $\tau > 0$ and $\tau \sinh a_\tau' = -1$
when $\tau < 0$. Then
\begin{equation}
s_\tau := \frac2{\pi}\,\int_0^{a_\tau}\frac{d\sigma}{
\sqrt{1 - \tau^2\cosh^2\sigma}}, \qquad \mbox{when} \quad \tau \in \Par^+,
\label{eq:1.3}
\end{equation}
and
\begin{equation}
s_\tau := \frac2{\pi}\,\int_0^{a_\tau'}\frac{d\sigma}{
\sqrt{1 -\tau^2\sinh^2 \sigma}}, \qquad \mbox{when}\quad \tau \in \Par^-.
\label{eq:1.5}
\end{equation}
We also use the following equivalent forms of these representations
\begin{equation}
s_\tau := \frac2{\pi}\,\int_0^{b_\tau}\frac{dx}{
\sqrt{\cos^2 x - \tau^2}}, \qquad \mbox{when} \quad \tau \in \Par^+,
\label{eq:1.31}
\end{equation}
where $b_\tau = \arccos \tau \in [0,\pi/2)$, $0 < \tau \leq 1$, and
\begin{equation}
s_\tau := \frac2{\pi}\,\int_0^{\pi/2}\frac{dx}{
\sqrt{\cos^2 x + \tau^2}}, \qquad \mbox{when}\quad \tau \in \Par^-.
\label{eq:1.51}
\end{equation}
These are effected by the changes of variables $\tau \cosh \sigma = 
\cos x$ and $\tau \sinh \sigma = \cos x$, respectively. 

\begin{proposition}
As a function of $\tau$, $s_\tau$ is monotone increasing when
$-\infty < \tau < 0$ and monotone decreasing when $0 < \tau \leq 1$.
\label{pr:3.1}
\end{proposition}
\begin{proof} It is obvious from (\ref{eq:1.51})
that $ds_\tau/d\tau > 0$ when $\tau < 0$. The fact that $s_\tau$ decreases 
when $\tau$ 
increases from $0$ to $1$ is less obvious and somewhat more
difficult to obtain; since we do not need it here we omit
the proof. 
\end{proof}

\medskip

We next consider the asymptotics of $s_\tau$ as $\tau$ approaches
various possible limiting values in $\Par$. The first two cases
are less important, but we record them anyway; namely, we have
\begin{equation}
s_\tau = - \frac{1}{\pi} \log \tau^2  + {\mathcal O} (1)
\qquad \mbox{as}\quad \tau \to 0
\label{eq:lims0}
\end{equation}
and also
\[
\lim_{\tau \nearrow 1} s_\tau = 1.
\]
On the other hand, the behaviour of $s_\tau$ as $\tau \to -\infty$
is fundamental to our analysis. To study this we introduce
an auxiliary function
\[
\gamma_\tau(t):= \tau\,\sigma_\tau \left(s_\tau\,t\right), \qquad \tau < 0;
\]
by design, $\gamma$ has period $2\pi$. 
\begin{lemma}
As $\tau \to -\infty$, $\gamma_\tau $ converges uniformly to the 
function $\cos t$.
More precisely, 
\begin{equation}
\gamma_\tau (t) = \cos t  + {\mathcal O}(|\tau|^{-2}) \quad  \mbox{and} \quad 
\del_\tau \gamma_\tau (t) = {\mathcal O}(|\tau|^{-3}),
\label{eq:3.1}
\end{equation}
as $\tau \to -\infty$. Simultaneously, in the same limit
\begin{equation}
s_\tau  = - \frac{1}{\tau} + {\mathcal O}(|\tau|^{-3}), \quad \mbox{and} \quad \del_\tau s_\tau  =  \frac{1}{\tau^2} + {\mathcal O}(|\tau|^{-4}).
\label{eq:3.2}
\end{equation}
\label{le:3.1}
\end{lemma}
\begin{remark}
The proof below easily extends to give full asymptotic expansions
for $s_\tau$ and $\gamma$ in powers of $1/\tau$ as $\tau \to -\infty$.
However, we only require the terms given in this statement.
\end{remark}
\begin{proof}  Set $\sigma_0 :=  \sigma(0)$. Since $|\sigma|
\leq -\sigma_0$, we may define the function $w(s)$ by 
\[
\sigma (s) = \sigma (0) \, \cos (\tau \, w(s)).
\]
To normalize it, we require that $w(0)=0$ and $w'(0) > 0$. 
The equation for $\sigma$, (\ref{eq:1.21}), becomes
\[
(\del_s w )^2  + 
\frac{\sinh^2 \sigma - \sinh^2 \sigma_0}{\sigma^2 -\sigma_0^2} = 0,
\] 
or equivalently
\begin{equation}
(\del_s w )^2   = \frac{\Phi(\sigma^2) - 
\Phi (\sigma_0^2)}{\sigma^2 -\sigma_0^2}
\qquad \mbox{where} \quad \Phi (t) : = \sinh^{\,2} \sqrt t.
\label{eq:eqg}
\end{equation}
(\ref{eq:1.21}) also implies that $|\sigma| \leq |\sinh \sigma| 
\leq |\tau|^{-1}$. Using the first term of the Taylor expansion
of $\Phi$ gives $\del_s w  = 1 + {\mathcal O}(|\tau|^{-2})$, and
so
\begin{equation}
w(s) = s \, (1+ {\mathcal O}(|\tau|^{-2})).
\label{eq:333}
\end{equation}
Next, the periodicity of $\sigma$ translates to the equality
\[
-2\pi = \tau w(2\pi s_\tau) = \tau(2\pi s_\tau  + {\mathcal O}(s_\tau 
|\tau|^{-2})),
\]
or equivalently,
\[
s_\tau = -\frac{1}{\tau} + {\mathcal O}(|\tau|^{-3}).
\]
Inserting this into the expression for $\gamma(t)$, and using
$\sigma_0 = {\mathrm{arcsinh}}\,(1/\tau) = 1/\tau + {\mathcal O}
(1/|\tau|^{-3})$ gives
\[
\gamma(t) = \cos t + {\mathcal O}(|\tau|^{-2}).
\]
This procedure may be continued, using this expansion for $\gamma$
(and hence for $\sigma$) to determine the next term in the expansion
for $s_\tau$, then using this in turn to get the next term
in the expansion for $\gamma$, and so on.  
\end{proof}

\section{Spectral analysis of the Jacobi operator}

The primary tool in this paper is a detailed analysis of the spectrum 
of the linearized mean curvature operator $\calL_\tau$ (which is
usually called the Jacobi operator) on the family of Delaunay surfaces 
$D_\tau$. We begin by describing the specific form of this operator and 
then review how Bloch wave theory may be applied in this context. This
theory reduces the spectral analysis to that of a countable collection of
continuous families of operators with discrete spectrum and it shows
that the spectrum of $\calL_\tau$ is a union of `bands' of absolutely
continuous spectrum. Geometric 
considerations allow us to identify specific solutions of $\calL_\tau
u = 0$, and this, in turn, allows us to track the location of 
some of these bands as $\tau \to -\infty$. 

\subsection{Coordinate expression}
It is well known that the Jacobi operator has the form
\[
\calL_\tau : = - \Delta_\tau - |A_\tau|^2,
\]
where $\Delta_\tau$ and $A_\tau$ are the Laplace operator and second
fundamental form on $D_\tau $, respectively. In terms of the 
isothermal parametrization (\ref{eq:1.1}), 
\begin{equation}
\calL_\tau = -\frac1{\tau^2\,e^{2\sigma}} \left( \del_{s}^2 + 
\del_{\theta}^2 - \tau^2 \cosh 2\sigma \right). 
\label{eq:4.1}
\end{equation}

We shall be using the periodicity of the coefficients of this operator
in a central way. However, it is more convenient to work with a family
of operators for which the period is fixed as $\tau$ varies.
Accordingly, we change variables, setting 
\[
s := s_\tau\,t,
\]
so that the coefficients of $L_\tau$ with respect to $(t,\theta)$ 
have period $2\, \pi$ for all $\tau \in \Par$. In these new coordinates, 
\[
\calL_\tau = -\frac1{s_\tau^2 \tau^2 e^{2\sigma}}
\left( \del_t^2 + s_\tau^2 \del_\theta^2 - 
s_\tau^2 \tau^2 \cosh(2\sigma) \right).
\]

Removing the factor $1/ (s_\tau^2 \tau^2 e^{2\sigma})$, we define
\[
L_\tau : = - \del_t^2 - s_\tau^2 \del_\theta^2 - s_\tau^2 \tau^2 
\cosh(2\sigma).
\]
Since we are only interested in some aspects of the spectral analysis 
of ${\calL}_\tau$, such as the number of negative eigenvalues, 
the existence of a nullspace, etc., we may concentrate on the study 
of the slightly simpler operator $L_\tau$.

\subsection{Spectral decomposition}
The operator $L_\tau$ has many symmetries, and these may be used
to reduce it to operators for which the spectral analysis
is more tractable. 

The first and most obvious reduction uses the rotational 
invariance in $\theta$. Thus, if $u(t,\theta) \in L^2(\R \times S^1)$,
we have the decomposition
\[
u(t,\theta) = \sum_{j\in {\mathbb Z}} u_j(t)e^{ij\theta},
\]
where each of the coefficients $u_j(t)$ is in $L^2(\R)$ and
$\sum ||u_j||^2 < \infty$. The operator $L_\tau$ induces the operator
\begin{equation}
L_{\tau,j} = - \del_t^2 + s_\tau^2 j^2 - s_\tau^2 \, \tau^2 \, \cosh (2\sigma) ,
\label{eq:5.2}
\end{equation}
on the $j^{{\mathrm{th}}}$ eigenspace. In fancier language,
the operator $L_\tau$ is reduced by the splitting
\[
L^2(\R \times S^1) = \bigoplus_{j\in {\mathbb Z}} L^2_j(\R)
\]
from this eigenspace decomposition into a direct sum of (self adjoint)
ordinary differential operators $L_{\tau,j}$. Note that $L_{\tau,j}
= L_{\tau,-j}$, even though the spaces $L^2_j$ and $L^2_{-j}$
are different (one corresponds to the eigenfunction $e^{ij\theta}$
and the other to the eigenfunction $e^{-ij\theta}$). It is clear that
\[
\mbox{spec}\,(L_\tau) = \bigcup_{j \in {\mathbb N}} \mbox{spec}\,(L_{\tau,j}).
\]
Furthermore, noting that 
\begin{equation}
L_{\tau,j} = L_{\tau,0} + s_\tau^2 j^2,
\label{eq:shift}
\end{equation}
we may as well restrict attention to the operator $L_{\tau,0}$ and
the single Hilbert space $L^2(\R)$. 

To analyze this last operator we use the technique of Bloch waves,
also known as Flocquet theory. We refer to \cite{Ree-Sim} and \cite{Mag}
for details, cf.\ also \cite{Maz-Pol-Uhl}. This relies on a direct integral 
decomposition of $L^2(\R)$ defined using the Fourier-Laplace transform. 
Given any function $f$ in the Schwartz space ${\mathcal S}$ on $\R$, 
define
\[
\hat{f}(t,\alpha) = \sum_{j \in {\mathbb Z}} f(t+2\pi j) \, e^{- i\alpha j}.
\]
We may as well assume that $-\pi \leq \alpha \leq \pi$. This function 
satisfies $\hat{f}(t+2\pi,\alpha) = e^{i\alpha}\hat{f}(t,\alpha)$,
and so we lose no information by restricting $t$ to lie in $[-\pi,\pi]$.
It is straightforward to check that
\[
\|f\|_{L^2(\R)}^2 = 2 \, \pi \, \|\hat{f}\|_{L^2([-\pi,\pi]^2)}^2,
\]
which means that the map $f \mapsto \hat{f}$ may be extended
as an isometry from $L^2(\R)$ into $L^2([-\pi,\pi]^2)$. In fact it is
an isometry, since there is an inversion formula given by
\[
f(t) = \frac{1}{2\pi}\int_{-\pi}^\pi \hat{f}(\bar{t},\alpha) \, e^{i\alpha j} \, d\alpha,
\]
where $\bar{t} \in [-\pi,\pi]$ and $t = \bar{t} + 2\pi j$.
Altogether, this gives the direct integral decomposition
\[
L^2(\R) = \int^{\oplus}_{\alpha \in [-\pi,\pi]} L^2_\alpha([-\pi,\pi]);
\]
the Hilbert spaces $L^2_\alpha$ appearing here are all the same,
but the index reflects the fact that they arise as the closures of
the space of continuous (or $H^1$) functions $u$ which satisfy $u(\pi) = 
e^{i\alpha}u(-\pi)$. The reason for this distinction will be evident soon. 

The importance of this decomposition for us is that the operator $L_{\tau,0}$
reduces further into a direct integral of operators $H_{0,\alpha}(\tau)$
induced on each one of the spaces $L^2_\alpha$. More specifically,
$H_{0,\alpha}(\tau)$ has the same expression as $L_{\tau,0}$, but its
domain is the space of $H^2$ functions on $[-\pi,\pi]$ satisfying 
the quasiperiodic boundary conditions $u(\pi) = e^{i\alpha}u(-\pi)$, 
$u'(\pi) = e^{i\alpha}u'(-\pi)$. 

As an aside about notation, we will be considering the reduction of
the full operator $L_\tau$ to the subspace of $\alpha$-quasiperiodic
functions, and shall denote this operator by $H_\alpha(\tau)$.
Its restriction to the $j^{\mathrm{th}}$ eigenspace of the cross-sectional 
operator is denoted $H_{j,\alpha}(\tau)$, and since
\begin{equation}
H_{j,\alpha}(\tau) = H_{0,\alpha}(\tau) + s_\tau^2 j^2,
\label{eq:redj}
\end{equation}
it will suffice at first to consider only the case $j=0$. Notice
that the subscripts $j$ and $\alpha$ refer to the symmetries
of this reduced subspace (i.e.\ $j^{\mathrm{th}}$ eigenvalue on the
cross-section and quasiperiodicity $\alpha$); the parameter $\tau$
has been elevated from subscript status to ground level, since 
our ultimate concern is 
with the spectral flow of these operators as a function of $\tau$.

Continuing the main thread of the discussion, the spaces $L^2_\alpha$ are 
all the same, but the subscript $\alpha$ is meant to remind us 
that the domains of the restrictions of $L_{\tau,0}$ to each one of 
them is different. Since these are self-adjoint boundary
conditions, each $H_{0,\alpha}(\tau)$ has discrete spectrum:
\[
\mbox{spec}\,(H_{0,\alpha}(\tau)) = \{\lambda_0(\tau,\alpha) \leq 
\lambda_1(\tau,\alpha) \leq \ldots \}.
\]
The $k^{\mathrm{th}}$ eigenvalue map can be considered as a map:
\[
\lambda_k(\tau,\cdot): S^1 \longrightarrow \R;
\]
its image is an interval $B_k(\tau)$, which is called the 
$k^{\mathrm{th}}$ band; the union of all these intervals constitutes
the spectrum of $L_{\tau,0}$:
\[
\mbox{spec}\,(L_{\tau,0}) = \bigcup_{k \in {\mathbb N}} B_k(\tau).
\]

A bit more is known about these bands, and the `band functions'
$\lambda_k$. First, since $L_{\tau,0}$ is an ordinary differential 
operator, the space of solutions of $(L_{\tau,0} - \lambda)\phi = 0$
(with no growth restrictions) is precisely two-dimensional.
Suppose that $\lambda$ is in the spectrum, i.e.\ $\lambda = 
\lambda_k(\tau,\alpha)$ for some $k$ and $\alpha$. 
The corresponding eigenfunction satisfies $\varphi(t+2\pi) = 
e^{i\alpha}\varphi(t)$; since the coefficients of the operator
are real, $\overline{\phi}$ is in the same eigenspace and
satisfies $\overline{\phi}(t+2\pi) = e^{-i\alpha}\overline{\phi}(t)$.
This implies that $\lambda_k(\tau,-\alpha) = \lambda_k(\tau,\alpha)$,
and hence we may as well restrict $0 \leq \alpha \leq \pi$. 
In fact, it is known \cite{Mag} that the band functions with {\it even} index 
are nondecreasing while the band functions with 
{\it odd} index are nonincreasing in this interval, so that
\[
\lambda_0(\tau,0) \leq \lambda_0(\tau,\pi) \leq \lambda_1(\tau,\pi)
\leq \lambda_1(\tau,0) \leq \lambda_2(\tau,0) \ldots.
\]
This means that the bands are nonoverlapping, except perhaps
at their endpoints, and we have
\[
\begin{array}{rllll}
B_{2k}(\tau) & = & \ds [\lambda_{2k}(\tau,0),\lambda_{2k}(\tau,\pi)]
\\[3mm]
B_{2k+1}(\tau) & = & \ds [\lambda_{2k+1}(\tau,\pi), \lambda_{2k+1}
(\tau, 0)].
\end{array}
\]

No band $B_k(\tau)$ may reduce to a point for any
$\tau \in \Par$. The reason is simply that if this were the case,
then the band function $\lambda_k(\tau,\alpha)$ would be
constant in $\alpha$, say $\lambda_k(\tau,\alpha) \equiv 
\lambda_k^0(\tau)$,
and hence there would be an infinite dimensional space of solutions
to the equation $L_{\tau,0}u = \lambda_k^0(\tau)u$, which is at
odds with the fact that $L_{\tau,0}$ is an ordinary differential
operator.
 
For any $\lambda \in \R$, the space of solutions of $L_{\tau,0}u =
\lambda u$ is two-dimensional. When $\lambda \in \cup_{k \in {\mathbb N}}
B_k(\tau)$, then the eigenspace contains two quasi-periodic (and
hence bounded) functions, whereas if $\lambda$ lies outside the
spectrum, then the corresponding set of solutions has a 
basis $\{u_\pm\}$, where $u_+$ decays 
exponentially as $t \to -\infty$ and grows exponentially as $t \to \infty$, 
while $u_-$ grows exponentially as $t \to -\infty$ and 
decays exponentially as $t \to \infty$. 

The main task lying ahead is to understand the dependence on $\tau$ 
of these bands $B_k(\tau)$. While it is certainly not possible to 
do this explicitly, we will use the special geometric 
solutions of $L_\tau u = 0$, the existence and nature of which 
give some nontrivial information about the bands $B_k(\tau)$, 
$k = 0, 1, 2$. 

\subsection{Geometric Jacobi fields}
Solutions of $\calL_\tau u = 0$ (or equivalently, $L_\tau u = 0$)
are called Jacobi fields; they correspond to variations of the Delaunay 
surface $D_\tau$ which preserve the mean curvature to second
order. We are particularly interested in solutions which
defined on all of $D_\tau$ which grow at most polynomially.
In general, Jacobi fields may not correspond to actual CMC deformations,
but as we describe in a moment, these correspond to global CMC
deformations which arise from explicit geometric motions.

The Bloch wave theory of the previous subsection, cf.\ \cite{Maz-Pol-Uhl}, 
\cite{Kus-Maz-Pol} can be used to deduce the Fredholm
properties and asymptotic expansions as $t \to \pm \infty$ 
for solutions of $L_\tau u = 0$. A consequence of this development
is the fact that Jacobi fields have definite exponential rates
of growth or decay. In particular, one proves that there is
at most a finite dimensional family of Jacobi fields $u$
which satisfy the bounds $|u|e^{-\epsilon |t|} \leq C$ for
some (and in fact any) sufficiently small $\epsilon$. Somewhat
remarkably, it is possible to identify every one of these
temperate Jacobi fields because all of them correspond to explicit
CMC deformations of $D_\tau$.

There are three types of global CMC deformations: those
corresponding to translations of $D_\tau$, which are a 
three-dimensional family, those corresponding to rotations
of $D_\tau$, which are a two-dimensional family, and a
final one-dimensional family corresponding to changing
the Delaunay parameter $\tau$. The infinitesimal variations
corresponding to each of these span the six-dimensional
space of temperate Jacobi fields. We now describe the
distinguished three-dimensional subspace of these which
are bounded. 

We temporarily revert to the isothermal $(s,\theta)$ coordinate
system.

\begin{itemize}
\item There is a Jacobi field $\Phi_\tau^{0}(s)$ which corresponds
to translating $D_\tau$ along its axis. It is obtained by projecting 
the constant Killing field $(0,0,1)$ generating this translation
over the normal vector field to $D_\tau$. This yields
\[
\Phi_\tau^{0} : =  \del_s \sigma_\tau.
\]
This Jacobi field is rotationally invariant and is also
periodic in $s$, hence bounded. In addition, $\Phi^{0}_\tau$ has 
precisely two nodal domains on the portion of $D_\tau$ where 
$s \in [0, 2 \, \pi \, s_\tau]$. 
\item The two Jacobi fields corresponding to translating $D_\tau$ 
orthogonally to its axis are given by projecting the constant Killing 
fields $(1,0,0)$ and $(0,1,0)$ on the unit normal field. These functions 
are again periodic in $s$; appropriate (complex) linear 
combinations of them have the form $\Phi_\tau^{\pm 1}(s)e^{\pm i \theta}$.
By closer inspection, we actually see that
\[
\Phi_\tau^{\pm 1} = \cosh \sigma \, e^{\pm i\theta}
\qquad \mbox{when}\quad \tau \in \Par^+,
\]
and
\[
\Phi_\tau^{\pm 1} =  \sinh \sigma \, e^{\pm i\theta}
\qquad \mbox{when} \quad \tau \in \Par^-.
\]
Notice that $\cosh \sigma$ does not change sign  while
$\sinh \sigma$ has precisely two nodal domains when $s \in
[0, 2 \, \pi \, s_\tau]$. 
\end{itemize}
\begin{remark}
There are three other geometric Jacobi fields, corresponding
to rotating $D_\tau$ around its axis and varying the Delaunay
parameter. These are not periodic but instead
grow linearly, hence do not enter our considerations below.
Therefore we omit further discussion of them.
\end{remark}

\subsection{The spectral flow}
We shall henceforth always reduce to $\alpha$-quasiperiodic
functions, and so it is no longer necessary to think in terms 
of the operator $L_\tau$ on the complete Delaunay surface $D_\tau$,
nor the band structure of its spectrum;  rather, we focus
on the operators $H_\alpha(\tau)$, which are simpler to analyze
because they have discrete spectrum. The ultimate goal now is
to determine when these operators become unstable, i.e.\ when 
they have nonzero index (number of negative eigenvalues). 
We know this does not occur when $\tau > 0$, nor when 
$\tau_* < \tau < 0$. In this section we analyze the spectral
flow of $H_\alpha(\tau)$ as $\tau$ decreases to $-\infty$. 

The eigenvalues of this operator take the form
$\lambda_{kj}(\tau) = 
\lambda_k(\tau,\alpha) + s_\tau^2 j^2$, $j,k = 0, 1, 2, \ldots$,
and so this study divides into essentially two parts. First
we obtain some information about the values of $\lambda_k(\tau,\alpha)$
when $k = 0,1,2$; following this
we turn to the problem of when $\lambda_{kj}(\tau) < 0$. 

We take up the first of these tasks in the next two subsections.
Because $\lambda_k(\tau,\alpha)$ always lies between $\lambda_k(\tau,0)$
and $\lambda_k(\tau,\pi)$, we focus most of the attention on
functions satisfying periodic or antiperiodic boundary conditions.

\subsubsection{The spectrum of $H_{0,\alpha}(\tau)$ when $\tau > 0$}
Although not needed later, we first show that there is no spectral flow
for the operators $H_\alpha (\tau ) $ as $\tau$ increases from $0$ to $1$. 

\begin{proposition}
Suppose that $\tau \in \Par^+$. Then the bottom of first band
$B_0(\tau)$ occurs at $\lambda = -s_\tau^2$, i.e. 
\[
\lambda_0(\tau,0) =  -s_\tau^2. 
\]
In addition, $0$ is either at the top of the second band
or bottom of the third, i.e.
\[
{\mathrm{either}} \qquad 
\lambda_{1}(\tau,0) = 0 \qquad {\mathrm{or}}
\qquad \lambda_{2}(\tau, 0) = 0.
\]
\label{pr:5.1}
\end{proposition}
\begin{proof}
We have indicated that $L_\tau(\cosh \sigma e^{i\theta}) = 0$, 
and it satisfies periodic boundary conditions, so that 
$H_{1,0}(\tau)(\cosh \sigma) = 0$. By (\ref{eq:redj}), 
this is the same as
\[
H_{0,0}(\tau)(\cosh \sigma) = -s_\tau^2 \cosh \sigma.
\]
But since $\cosh \sigma$ is everywhere positive, it must correspond
to the ground state eigenvalue.

Next, the function $\Phi_\tau^0 = \del_s \sigma $ is a solution
of $L_{\tau} w = 0$, and again satisfies periodic boundary conditions;
moreover, it has precisely two nodal regions in a period domain, and so 
it must correspond either to the second or third eigenvalue of 
$H_{0,0}(\tau)$, i.e.\ either to $\lambda_1(\tau,0)$ or $\lambda_2(\tau,0)$, 
as stated. Unfortunately, it is not clear from the evidence at hand
which of these is actually the case.  
\end{proof}

\begin{corollary} When $\tau > 0$, the bands $B_k(\tau)$, $k \geq 3$, 
are strictly contained in the positive half-line $\R^+$, while
$B_2(\tau) \subset \overline{\R^+}$ (possibly in the 
open half-line). In any event, the operator $H_{0,\alpha}(\tau)$ has no 
spectral flow as $\tau$ varies in $\Par^+$.
\end{corollary}
The first statement follows at once from the previous Proposition.
Obviously, since $\lambda_k(\tau, \alpha) \geq -s^2_\tau$, we have 
$\lambda_{kj}(\tau)\geq 0$ for all $j \geq 1$, $k \geq 0$ and all $\tau 
\in [0, \pi]$. This proves the last statement.

\subsubsection{The spectrum of $H_{0,\alpha}(\tau)$ when $\tau < 0$}
In contrast to the situation for unduloids, there is more to say about 
the spectrum of $H_\alpha (\tau)$ for nodoids. 
\begin{proposition}
Let $\tau \in \Par^-$. Then for $\alpha \in [0, \pi]$ we have
\[
-2s_\tau^2 + \left(\frac{\alpha}{2\pi}\right)^2 - \tau^2 s_\tau^2
\leq \lambda_0 (\tau, \alpha)
\leq \left(\frac{\alpha}{2\pi}\right)^2 - \tau^2 s_\tau^2;
\]
in addition,
\[
\lambda_1 (\tau, 0) =  - s_\tau^2, \qquad \lambda_2 (\tau, 0) = 0.
\]
\label{pr:5.2}
\end{proposition}
\begin{proof} First consider the identities for $\lambda_1$ and 
$\lambda_2$.  We already know that $L_{\tau}(\sinh \sigma \, e^{i\theta}) 
= 0$ and this function again satisfies periodic boundary conditions,
so that $H_{0,0}(\tau)(\sinh \sigma) = -s_\tau^2 \sinh \sigma$.
This function has two nodal regions, hence $-s_\tau^2$ must equal
either the second or third eigenvalue of $H_{0,0}(\tau)$. On the other hand,
we also know that $L_{\tau}\del_s \sigma = 0$, and this function is
again periodic, so $0$ is also in the spectrum of $H_{0,0}(\tau)$. 
Since $\del_s \sigma$ again has only two nodal regions, $0$ must also 
be either the second or the third eigenvalue. Putting these statements 
together shows that $-s_\tau^2$ is the second eigenvalue and $0$ is the 
third eigenvalue of $H_{0,0}(\tau)$.

\medskip

To obtain the bounds on the bottom eigenvalue, note that 
the potential in $L_{\tau,0}$ satisfies the estimates
\[
-s_\tau^2 (\tau^2 + 2) \leq -s_\tau^2 \tau^2 \, \cosh (2 \sigma) 
\leq  -s_\tau^2 \tau^2 
\]
since  
\[
\tau^2 \, \cosh (2\sigma) = \tau^2 +2 \, \tau^2 \sinh^2 \sigma, 
\]
and $\tau^2 \sinh^2 \sigma \leq 1$. The estimate for $\lambda_0(\tau,\alpha)$
is then straightforward by monotonicity, but cf.\ below for the precise 
form of the eigendata of these $\alpha$-quasiperiodic problems when the 
potential is replaced by a constant.  \end{proof}

\medskip

We have now shown that 
\[
B_0(\tau) \cup B_1(\tau) = [\lambda_0(\tau,0),\lambda_0(\tau,\pi)] 
\cup [\lambda_1(\tau,\pi),-s_\tau^2] \subset (-\infty,0),
\]
and
\[
B_2(\tau) = [0,\lambda_2(\tau,\pi)].
\]
Recalling our earlier remark that no band reduces to a point, 
we have $\lambda_2(\tau,\pi) > 0$ and hence 
$B_k(\tau) \subset \R^+$ when $k \geq 3$. 

\begin{remark}
The main conclusion of this discussion is that only the bands $B_0(\tau)$ 
and $B_1(\tau)$ lie in the negative half-line; all other bands 
are wholly contained in the positive half-line.
\end{remark}

It will also be necessary later to have more refined information
about the behaviour of the $\alpha$-quasiperiodic eigenvalues
as $\tau \to -\infty$. The key observation is that $L_{\tau,0}$ 
converges uniformly in this limit to $L := - \del_t^2 - 1$, and 
hence the spectrum of $L_{\tau, 0}$ converges to that of $L$. 
Clearly $\mbox{spec}\,(L) = [-1,\infty)$, but although $L$ has 
constant coefficients we may still perform the Bloch wave analysis to 
decompose this spectral ray into an infinite union of spectral bands. 
In fact, the solutions of $Lu = \lambda u$ with $u(2\pi)= e^{i\alpha}u(0)$, 
$u'(2\pi)= e^{i\alpha}u'(0)$ are given by $e^{i((\alpha/2\pi) \pm  k)t}$, 
and so for $\alpha \in [0, \pi]$, we have  
\[
\lambda_k(\alpha) = \left\{ 
\begin{aligned} \left(\frac{\alpha}{2\pi} +k \right)^2 - 1
\qquad \mbox{when}\quad k \quad \mbox{is even} \\[3mm]
\left( \frac{\alpha}{2\pi}  -k\right)^2 - 1
\qquad \mbox{when}\quad k \quad \mbox{is odd},
\end{aligned} \right.
\]
and the corresponding eigenfunctions are given by 
\[
\varphi_k(\alpha) = \left\{ 
\begin{aligned} \ds e^{i (\frac{\alpha}{2\pi} -k ) t}
\qquad \mbox{when}\quad k \quad \mbox{is even} \\[3mm]
\ds e^{i ( \frac{\alpha}{2\pi}+k ) t}
\qquad \mbox{when}\quad k \quad \mbox{is odd}.
\end{aligned} \right.
\]
This gives the band structure
\[
B_0 = [-1, -3/4], \quad B_1 = [-3/4, 0], \quad B_2 = [0, 5/4], \qquad 
\mbox{etc.}
\]
The apparent lack of smoothness in the band functions $\lambda_k(\alpha)$
at $\alpha = 0, \pi$ is due to the absence of gaps between these
bands. 

We now use Lemma~\ref{le:3.1} to perturb off this limiting situation.

\begin{proposition}
For $\alpha \in [0, \pi]$, let $(\varphi_k(\alpha),\lambda_k(\alpha))$ and 
$(\varphi_k(\tau,\alpha),\lambda_k(\tau,\alpha))$ denote the
eigenfunctions and eigenvalues for $L_\tau$ and $L_{\tau,0}$, respectively.
Then for $\tau$ sufficiently negative, 
\begin{equation}
\lambda_k (\tau,\alpha) = \lambda_k(\alpha) + {\mathcal O}(|\tau|^{-2}), 
\quad \quad \del_\tau \lambda_k (\tau,\alpha) =  {\mathcal O}(|\tau|^{-3}),
\label{eq:pertl}
\end{equation}
and moreover, 
\begin{equation}
\varphi_k(\tau,\alpha) = \varphi_k(\alpha) + {\mathcal O}(|\tau|^{-2})
\label{eq:perte}
\end{equation}
uniformly along with all derivatives on $0 \leq t \leq 2\pi$. 
\label{pr:5.3}
\end{proposition}
\begin{proof} The fact that the limits of $\lambda_k(\tau,\alpha)$
and $\varphi_k(\tau,\alpha)$ converge to $\lambda_k(\alpha)$ and
$\varphi_k(\alpha)$ is clear from general theory. To get the estimate
on the eigenvalue, use Lemma~\ref{le:3.1} to get that 
\begin{equation}
\begin{array}{rcl}
s_\tau^2 \, \tau^2 \cosh(2\sigma(s_\tau t))  & = &  1+
{\mathcal O}(|\tau|^{-2}), \qquad \mbox{and} \\
\del_\tau \left( s_\tau^2 \, \tau^2 \cosh(2\sigma(s_\tau t)) \right) & = & 
{\mathcal O}(|\tau|^{-3}).
\end{array}
\label{eq:pot}
\end{equation}
Assuming that the functions $\varphi_k(\tau,\alpha)$ are
normalized to have $L^2$ norm equal to $1$, we use the standard
formula from eigenvalue perturbation theory
\[
\del_\tau \lambda_k(\tau,\alpha) = \big\langle \del_\tau L_{\tau,0}\,  
\varphi_k(\tau,\alpha),\varphi_k(\tau,\alpha) \big\rangle_{L^2}
=  {\mathcal O}(|\tau|^{-3}),
\]
which follows directly from (\ref{eq:pot}). Integrating from 
$-\infty$ to $\tau$ gives the first part of (\ref{eq:pertl}); then a standard
perturbation argument yields (\ref{eq:perte}).
\end{proof}

\subsubsection{Spectral flow of $H_\alpha(\tau)$} 
We now let $H_\alpha(\tau)$ denote the operator $L_\tau$ 
acting on $\alpha$-quasiperiodic functions. (Thus it is 
rotational invariant, and its reductions to the eigenspaces
of the cross-section are the operators $H_{j,\alpha}(\tau)$.)
The motivation for this whole paper is the fundamental observation
that $H_\alpha(\tau)$ has a nontrivial spectral flow as $\tau$
decreases from $0$ to $-\infty$. To see this, recall that we have 
already determined that its eigenvalues are of the form 
\[
\lambda_k(\tau,\alpha)+ s_\tau^2 j^2, \qquad j,k = 0, 1, 2, \ldots 
\]
In addition, when $\alpha\in [0, \pi]$, $\lambda_0(\tau,\alpha)$
stays strictly negative and uniformly bounded away from zero as $\tau 
\to -\infty$, and the same is true for $\lambda_1(\tau,\alpha)$ 
for $\alpha\in (0, \pi]$ (but note that $\lambda_1(\tau,0) = -s_\tau^2
\to 0$). In particular, we see that $\lambda_0(\tau,\alpha) + s_\tau^2 j^2$
is positive for all $j$ (because of (\ref{eq:lims0})), but eventually
becomes negative for each fixed $j$. Thus more and more eigenvalues, 
corresponding to higher and higher eigenmodes on the cross-section, 
change sign from positive to negative. 

The following two propositions make this more precise. The first
is valid for any $j \geq 2$ and for `intermediate' values of $\tau$.
The second gives a much more accurate estimate, but is only valid
when $|\tau|$ is large.

\begin{proposition}
Let $j \geq 2$. If the interval 
$B_0(\tau) + j^2 s_\tau^2 $ contains $0$, then $B_0(\tau) + 4j^2 s_\tau^2$
is entirely contained in the positive axis $(0,\infty)$.
\label{pr:ifif}
\end{proposition}
\begin{proof} It follows from Proposition~\ref{pr:5.2} that 
\[
\lambda_0 (\tau, \alpha )\geq - (\tau^2+2) \, s_\tau^2 .
\]
Hence, $B_0(\tau) + j^2 s_\tau^2 $ cannot contain $0$
 if $\tau^2 < 2$ and $j \geq 2$. Therefore, we can assume that 
$\tau < - \sqrt 2$.

Now, suppose the result were to fail. Then for some $j$
we would have both
\[ 
\lambda_0(\tau,\pi) + j^2 s_\tau^2 \geq 0 \qquad \mbox{and}\qquad
\lambda_0(\tau,0) + 4j^2 s_\tau^2 < 0.
\]
This gives
\[
\lambda_0(\tau,0) - 4\lambda_0(\tau,\pi) \leq 0.
\]
On the other hand, from Proposition~\ref{pr:5.2}, 
\[
-(\tau^2 + 2)s_\tau^2 \leq \lambda_0(\tau, 0) < \lambda_0(\tau,\pi) 
\leq \frac14 - \tau^2 s_\tau^2,
\]
which implies that
\[
\lambda_0(\tau,0) - 4\lambda_0(\tau,\pi) \geq -1 + 3\tau^2 s_\tau^2 - 
2s_\tau^2.
\]
Combining these two inequalities we get
\[
s_\tau^2 \leq \frac{1}{3\tau^2 - 2}.
\]
But from (\ref{eq:1.51}), 
\[
s_\tau^2 \geq \frac{1}{1+\tau^2}.
\]
These two inequalities are incompatible when $\tau < -\sqrt{2}$. 
\end{proof}

\subsection{The index of $L_\tau$}
We now apply the information we have about the spectrum of $H_\alpha(\tau)$
to the spectral flow and index of $L_\tau$ restricted to 
to the subspace functions which are $T_{j,\alpha}$-symmetric and also
invariant under the reflection $(t,\theta) \mapsto (-t,-\theta)$.
Equivalently, we let $L_\tau$ act on functions $u$ which satisfy
\[
u(t + 2 \pi ,\theta ) = u (t , \theta +\alpha), \quad  
u(t , \theta + 2 \pi/j) = u (t , \theta) \quad \mbox{and} \quad 
u(-t,-\theta) = u(t, \theta)
\]
for all $(t,\theta) \in [-\pi,\pi]\times S^1$. As before,
we denote the $H^2$ closure of this space by $H^2_{j,\alpha}$. 

Observe that the Fourier decomposition of any 
function $u \in H^2_{j,\alpha}$ reads
\[
u(t, \theta) = \sum_{n\in {\mathbb Z}-\{0\}} 
u_n (t) \, e^{i \, n \, j \, \theta}
\]
where $u_n$ satisfy $u_n (t+ 2 \pi) = 
e^{nj \alpha} \, u_n (t)$ and $u_n(-t)= -u_{-n} (t)$.
 
It is easy to see that the spectrum of $L_\tau$ acting 
on the space above is given by 
\[
\mbox{spec} (L_\tau) = \bigcup_{k\geq 0} \, \bigcup_{n\geq 1} 
\left( s_\tau^2 \, n^2\, j^2 + B_k (\tau)\right).
\]

\medskip

Now, for $j \geq 2$, $\tau \in \Par^-$ and $|\alpha| \leq \pi/j$, 
define $I_{j,\alpha}(\tau)$ to be the number of negative eigenvalues 
of the operator $L_\tau$ acting on $H^2_{j,\alpha}$. 
We have established that there exists a $\tau_* < 0$ such that
$I_{j,\alpha}(\tau) = 0$ when $\tau \in (\tau_*,0)$. As $\tau$
decreases, $\tau_*$ is the first value at which an eigenvalue of
$H(\tau)$ crosses zero; since we are restricting to some proper
subspace of functions, the spectrum of $H_{j,\alpha}(\tau)$
may remain positive on a larger interval. Thus define 
\[
\tau_{j,\alpha} = \inf\,\{\tau\,:\,I_{j,\alpha}(\tau') = 0 \, : \, 
\tau' \in (\tau, 0) \};
\]
hence all eigenvalues of $H_{j,\alpha}(\tau)$ remain positive when
$\tau \in (\tau_{j,\alpha},0)$. 
So, for all $\tau \in (\tau_{j, \alpha})$, the index of $L_\tau$ is 
$0$ and for some $\tau$ which is slightly smaller than 
$\tau_{j, \alpha}$, the index of $L_\tau$ is 
at least $1$. Observe that $\tau_{j,\alpha}$ is a root of the equation 
$\lambda_0(\tau, j \, \alpha) + s_\tau^2 j^2 = 0$. The 
next result shows that we have a good control on 
$\tau_{j, \alpha}$ when $j$ is large enough.
\begin{proposition}
There exists a $j_0 \geq 2$ such that if $j \geq j_0$, and if
$|\alpha| \leq \pi$, then
\[
\tau_{j,\alpha} =- \frac{j}{\sqrt{1-(\alpha/2\pi)^2}} + 
{\mathcal O}(j^{-1}).
\]
Moreover $I_{j,\alpha}(\tau)=0$ for $\tau > \tau_{j,\alpha}$ 
and $I_{j,\alpha}(\tau)\geq 1$  for all $\tau < \tau_{j,\alpha}$. 
\label{pr:6.2}
\end{proposition}
\begin{proof} The number $\tau_{j,\alpha}$ corresponds to
the value of $\tau$ for which $\lambda_0(\tau, j \, \alpha) + s_\tau^2 
j^2 = 0$. The estimate here follows from (\ref{eq:pertl}) and the explicit
expression for $\lambda_k(\alpha)$. The point of requiring $j$ to
be large is that this forces $\tau$ also to be large, and then
we can use these asymptotics results. 
\end{proof}

\medskip

For general $j \geq 2$ our control on $\tau_{j, \alpha}$ is weaker. 
The proposition below is immediate from Proposition~\ref{pr:5.2}
\begin{proposition} 
For $j \geq 2$, we have $\tau_{j,0},\tau_{j,0}' \in (-j,-\sqrt{j^2-2})$.
\label{pr:6.1}
\end{proposition}

The next result is a consequence of the proof of Lemma~\ref{le:3.1}, and is 
the basis for proving that when $j$ is large enough, the
eigenvalue crossings are transversal.
\begin{proposition} Fix $\alpha \in [-\pi,\pi]$. 
There exists $j_0\geq 2$ such that if $j \geq j_0$ 
and if $\lambda_k(\tau,\alpha) + s_\tau^2 j^2=0$, $k=0,1$, then
\[
\del_\tau (\lambda_k(\tau,\alpha) + s_\tau^2 j^2 ) < 0 .
\] 
\label{pr:trans}
\end{proposition}
\begin{proof} First, from Lemma~\ref{le:3.1}, we know  
that $\del_\tau \lambda_k(\tau,\alpha) =  {\mathcal O}(|\tau|^{-3})$.
Also, by (\ref{eq:3.2}), $s_\tau^2 = 1/\tau^2 + {\mathcal O}(|\tau|^{-4})$,
and $\del_\tau s_\tau^2 = -2/\tau^3 + {\mathcal O}(|\tau|^{-5})$.
 Hence
\[
\del_\tau (\lambda_k(\tau,\alpha) + s_\tau^2 j^2 ) =  - 2j^2 \, \tau^{-3}
+ {\mathcal O}(|\tau|^{-3}).
\]
Then observe that, since $\lambda_k (\tau,\alpha) 
+ s_\tau^2 j^2=0$ (recall that $\lambda_k$ is bounded since
we only need consider $k=0,1$), then $s_\tau \sim j^{-1}\sim |\tau |^{-1}$, 
so $\del_\tau (\lambda_k(\tau,\alpha) + s_\tau^2 j^2 )$ is certainly 
negative when $j$ is large enough.
\end{proof}

\section{CMC deformations of nodoids}

We now employ the preceding results about the Jacobi operator
to deduce the existence of CMC surfaces with $T_{j,\alpha}$ symmetry.

\subsection{The mean curvature operator}

For any $\tau \in \Par$, parametrize $D_\tau $ by 
\[
X_\tau(t,\theta) := \frac{1}{2}\, \left(\tau \, e^{\sigma_\tau
(s_\tau t)}\, \cos\theta,\tau\, e^{\sigma_\tau(s_\tau t)}\, \sin\theta, 
\kappa_\tau(s_\tau t)\right).
\]
The unit normal at $X_\tau (t,\theta)$ is then
\begin{equation}
N_\tau (t,\theta) := \big(\tau \, \sinh\sigma_\tau(s_\tau t)\,
\cos\theta,\tau \,\sinh\sigma_\tau(s_\tau t)\, \sin\theta, -\del_s 
\sigma_\tau (s_\tau t)\big).
\label{eq:6.1}
\end{equation}

If $w$ is any function which is $\calC^2$ small, then let $D_w$ 
the image of the map 
\[
X_w = X_\tau  + w\, N_\tau.
\]
Note that if $w$ is $R_j$ or $S_\alpha$ or 
$T_{j,\alpha}$ symmetric, then $D_w$ has the same symmetries.

\medskip

A rather complicated nonlinear elliptic equation determines
when $D_w$ has mean curvature $1$. We write it in abbreviated form as 
\begin{equation} 
L_\tau w  + Q_\tau (w)=0, 
\label{eq:6.2}
\end{equation} 
where $L_\tau$ is the multiple of the Jacobi operator which we have 
been studying, and $Q_\tau$ is a second order nonlinear 
differential operator which vanishes quadratically at $w=0$. More
precise information about the structure of $Q_\tau$ is given 
in \cite{Maz-Pac}. Note that in the $(t,\theta)$ coordinate system,
$Q_\tau$ is $2\, \pi$ periodic.

\medskip

Because this is a nonlinear problem, we shall use the function spaces
$\calC^{k,\beta}_{j,\alpha}(\R \times S^1)$, defined for any $j\geq 2$ 
and $\alpha \in [-\pi/j, \pi/j]$ to contain all $\calC^{k,\beta}$
functions which are $T_{j,\alpha}$-symmetric and invariant under
the reflection $(t,\theta) \mapsto (-t,-\theta)$:
\[
\C^{k,\beta}_{j,\alpha}(\R\times S^1)  =  \left\{u \in 
\calC^{k,\beta}(\R\times S^1): \, u(t+2\pi,\theta) = 
u(t,\theta +\alpha ), \right.
\]
\[
\left. u(t,\theta + 2\,\pi/j ) =  u(t,\theta)
\quad \mbox{and}\quad u(-t,-\theta) = u (t, \theta)\right\}.
\]

Clearly 
\[
L_\tau :\calC^{2,\beta}_{j,\alpha}(\R\times S^1)
\longrightarrow \calC^{0,\beta}_{j,\alpha} ({\R}\times S^1)
\]
and
\[
Q_\tau: \calC^{2,\beta}_{j,\alpha}(\R\times S^1)\longrightarrow 
\calC^{0,\beta}_{j,\alpha}(\R\times S^1)
\]
are smooth. 

\medskip

Consider the quotient of $\R^3$ by the screw motion $S_\alpha$ 
along the $z$-axis, with translation length $2\pi$; this is a $2$-plane 
bundle over $S^1$ with holonomy $\alpha$. $D_\tau/S_\alpha$ is a compact 
submanifold of this space, and we shall construct the surfaces bifurcating 
from $D_\tau$ as perturbations of $D_\tau/S_\alpha$. 

\subsection{Bifurcations}

We are now in a position to prove the existence of families of
(immersed) CMC surfaces which bifurcate off the Delaunay
surfaces $D_\tau$.

\medskip

The proof of Theorem~\ref{th:2.1} follows from the general 
bifurcation theorem of Smoller and Wasserman \cite{Smo-Was},
cf.\ also Theorem~13.10 in \cite{Smo}.  
To apply this result, we require only the fact that, by definition, 
the index $I_{j,\alpha}(\tau)$ is $0$ for any
$\tau > \tau_{j,\alpha}$ and, by Proposition~\ref{pr:ifif}, the 
index $I_{j,\alpha}(\tau)$ is exactly $1$ for some 
$\tau < \tau_{j,\alpha}$, but close to $\tau_{j,\alpha}$. 
This change of multiplicity ensures
the existence of a bifurcation for the nonlinear problem
(\ref{eq:6.2}) in the space $\calC^{2,\beta}_{j,\alpha}(\R\times S^1)$
(modulo the screw motion $S_\alpha$). When $\alpha = 0$, this produces 
nonrotationally invariant CMC surfaces which are periodic and 
$R_j$--symmetric. When $\alpha \neq 0$ this produces nonrotationally invariant 
CMC surfaces which are $T_{j,\alpha}$--symmetric.
By Proposition~\ref{pr:6.1} we get the estimate on the location
of this bifurcation point. 

\medskip

One defect of this general theorem is that we obtain no information
about whether this bifurcation gives a smooth connected branch of
solutions. For this we require nondegenerate crossing of a simple
eigenvalue. 
 However, in Proposition~\ref{pr:trans}, we have verified this hypothesis when 
$j$ and hence $\tau$ is sufficiently large. Therefore, we obtain
Theorem~\ref{th:2.2} from the theorem of Crandall and
Rabinowitz \cite{Cra-Rab}, cf.\ also Theorem~13.5 in \cite{Smo}. 

\medskip

We conclude by noting that there are many unresolved questions
concerning the surfaces we have produced. The most obvious one 
concerns the existence of a second bifurcation which arises when
$\lambda_1(\tau,\alpha) + s_\tau^2 j^2$ crosses zero. This should
not be difficult to obtain, and requires only a slight elaboration
of the techniques and estimates we have been using.
One complication here is how to separate the second bifurcation
for some smaller value of $j$ occurring at the same value
of $\tau$ as the first bifurcation for a larger value of $j$.
We have not pursued this because our information about these
surfaces is so limited, since we have proved their existence
using an abstract functional analytic technique. 
The most
interesting problem is to globalize this construction and find
a complete characterization of all immersed cylindrically bounded 
CMC surfaces with two ends.

\end{document}